\documentclass{commat}

\DeclareMathOperator{\height}{h}
\DeclareMathOperator{\ord}{ord}
\DeclareMathOperator{\SL}{SL}

\newcommand{\C}{{\mathbb{C}}}
\newcommand{\N}{{\mathbb N}}
\newcommand{\NN}{\mathcal{N}}
\newcommand{\OO}{\mathcal{O}}
\newcommand{\Q}{{\mathbb Q}}
\newcommand{\R}{{\mathbb{R}}}
\newcommand{\tild}{{\tilde{d}}}
\newcommand{\tilD}{{\widetilde{D}}}
\newcommand{\tilK}{{\widetilde{K}}}
\newcommand{\tils}{{\tilde{s}}}
\newcommand{\tilS}{{\widetilde{S}}}
\newcommand{\tilv}{{\tilde{v}}}
\newcommand{\Z}{{\mathbb Z}}

\title{%
An explicit bound of integral points on modular curves
}

\author{%
Yulin Cai
}

\affiliation{
    \address{%
    Institut de Math\'ematiques de Bordeaux, Universit\'e de Bordeaux 
	351, cours de la Lib\'eration 33405 Talence Cedex, France.
        }
    \email{%
    yulin.cai1990@gmail.com
        }
    }

\abstract{
In this paper, we give the constant $C$ in \cite[Theorem 1.2]{sha2014bounding} by using an explicit Baker's inequality, hence we obtain an explicit bound for the heights of the integral points on modular curves.
    }


    
\VOLUME{30}
\NUMBER{1}
\firstpage{161}
\DOI{https://doi.org/10.46298/cm.9389}

\begin{paper}

\section{Introduction}

Let~$X$ be a smooth, connected projective algebraic curve defined over a number field~$K$, and let ${x\in K(X)}$ be a non-constant rational function on~$X$. If~$S$ is a finite set of places of~$K$ (including all the infinite places), we call a point $P \in X(K)$ an \textit{$S$-integral point} if $x(P) \in \OO_S$, where ${\OO_S=\OO_{S,K}}$ is the ring of $S$-integers in~$K$. The set of $S$-integral points is denoted by $X(\OO_S,x)$.

According to the classical theorem of Siegel~\cite{siegle1929uber} the set $X(\OO_S,x)$ is finite if at least one of the following conditions is satisfied: 
\begin{align}
&\textrm{the genus } g(X)\geq 1;\\
\label{ethreepoles}
&\text{$x$ admits at least~$3$ poles in $X(\bar\Q)$}. 
\end{align}
Unfortunately, the existing proofs of this theorem for general curves are not effective, that is they do not imply any explicit expression bounding the heights of integral points. But for many pairs $(X,x)$, the effective proofs of this theorem were discovered by Baker's method, see \cite{bilu1995effctive}, \cite{bilu2002baker} and the references therein.

Sha \cite{sha2014bounding} considered the case where ${X=X_\Gamma}$ is the modular curve corresponding to a congruence subgroup~$\Gamma$ of $\SL_2(\Z)$, and ${x=j}$ is the $j$-invariant. 

To state his result, we introduce some notations. For a congruence subgroup $\Gamma$ as above, the number of cusps on $X_{\Gamma}$ is denoted by $v_\infty(\Gamma)$.
For a number field $K$, let $M_K$ be the set of all places of $K$, and $S \subseteq M_K$  a finite subset containing all infinite places. We put $d=[K:\mathbb{Q}]$ and $s = |S|$. Let $\mathcal{O}_K$ be the ring of integers of $K$. We define the following quantity
$$\Delta(N): = \sqrt{N^{dN}|D|^{\varphi(N)}}(\mbox{log}(N^{dN}|D|^{\varphi(N)}))^{d\varphi(N)}\times \left(\prod\limits_{\substack{v \in S\\  v\nmid \infty}}\mbox{log}\mathcal{N}_{K/\mathbb{Q}}(v)\right)^{\varphi(N)}$$
as a function of $N \in \N^+$, where $D$ is the absolute discriminant of $K$, $\varphi(N)$ is Euler's totient function, and the norm $\mathcal{N}_{K/\mathbb{Q}}(v)$ of a place $v$, by definition, is equal to $|\mathcal{O}_K/\mathfrak{p}_v|$ when $v$ is finite and $\mathfrak{p}_v$ is its corresponding prime ideal, and is set to be $1$ if $v$ is infinite.

Sha \cite{sha2014bounding} proved the following theorem.

\begin{theorem}[\cite{sha2014bounding} Theorem 1.2]
	\label{thsha}
	Let $\Gamma$ be of level $N$. If $v_\infty(\Gamma) \geq 3$, then
	\[
	    \height(j(P))\leq (CdsM^2)^{2sM}(\log(dM))^{3sM}\ell^{dM}\Delta(M)
	    \quad \textup{ for every } \ P \in X_\Gamma(\mathcal{O}_S,j),
	\]
	where $C$ is an absolute effective constant, $\ell$ is the maximal prime such that there exists $v\in S$ with $v|\ell$, or $\ell = 1$  if $S$ only contains infinite places,  and $M$ is defined as following: 
	$$M = \begin{cases}
	N & \text{if  $N$ is not a power of any prime;}\\
	3N & \text{if  $N$ is a power of $2$;}\\
	2N & \text{if $N$ is a power of an odd prime.}
	\end{cases}
	$$	
	\label{bound}
\end{theorem} 
(Here ${\height(\cdot)}$ is the standard absolute logarithmic  height defined on the set $\bar\Q$ of algebraic numbers.)

For certain applications it is useful to have an explicit value of the constant $C$ from Theorem~\ref{bound}. In this note we prove the following result.
\begin{theorem}
	The constant $C$ in Theorem~\ref{bound} can be taken to be $2^{14}$.
	\label{constant}
\end{theorem}
In the proof, we follow the main lines of Sha's argument, with some minor modifications. We calculate explicitly the implicit constants occurring therein.

For a number field $K$, $v \in M_K$, we define the valuation $|\cdot|_v$ on $K$ as following: for any $\alpha \in K$:
$$|\alpha|_v: = |\sigma(\alpha)|, \ \textrm{if $v$ is infinite with embedding $\sigma$};$$ 
$$|\alpha|_v: = \NN_{K/\Q}(v)^{-\ord_v(\alpha)/[K_v: \Q_v]}, \ \textrm{if $v$ is finite}.$$ 

\section{Upper bound of $S$-regulator} 
\label{S-regulator}

As before, for a number field $K$, and a finite subset $S \subseteq M_K$ containing all infinite places, we put $d=[K:\mathbb{Q}]$, $s = |S|$ and $r = s -1$. We fix $v_0 \in S$ and set
\[
    S' = S\setminus\{v_0\} = \{v_1, \dotsc, v_r\}.
\]
The $S$-regulator $R(S)$ is defined as
$$R(S) = |\det(d_{v_i}\log|\xi_k|_{v_i})_{1 \leq i,k\leq r}|,$$
where $d_{v_i} = [K_{v_i}: \Q_{v_i}]$ is the local degree of $v_i$ for each $i$, and $\{\xi_1,\cdots, \xi_r\}$ is a fundamental system of the $S$-units. It is independent of the choice of $v_0$ and of the fundamental system of $S$-units. We also denote by $\omega_K$ the number of roots of unity in $K$.

We set

\begin{equation*}
\zeta=\begin{cases}
\dfrac{(\log 6)^3}{2}&\mbox{if $d=2$,}\\
4\left(\dfrac{\log d}{\log\log d}\right)^3 &\mbox{if $d \geq 3$}.
\end{cases}
\end{equation*}
This $\zeta$ is better than the one in \cite[Proposition 4.1 and Corollary 4.2]{sha2014bounding}, and can make these results valid, see \cite[Theorem and Corollary 2]{voutier1996effective}.

\begin{lemma}
	We have
	$$0.1 \leq R(S)\leq h_KR_K\prod\limits_{\substack{v \in S \\ v  \nmid \infty}}\log\mathcal{N}_{K/\mathbb{Q}}(v),$$
	$$R(S)\leq\dfrac{\omega_K}{2}\left(\dfrac{2}{\pi}\right)^{r_2}\left(\dfrac{e\log |D|}{4(d-1)}\right)^{d-1}\sqrt{|D|}\prod\limits_{\substack{v \in S \\ v  \nmid \infty}}\log\mathcal{N}_{K/\mathbb{Q}}(v),$$
	where $e$ is the base of the natural logarithm, $r_2$ is the number of complex embeddings of $K$, and $D$ is the absolute discriminant of $K$.
	\label{bound RS}
\end{lemma}
\begin{proof}
	For the first inequality see  \cite[Lemma 3]{bugeaud1996bounds}. One may remark that the lower bound $R(S)\geq 0.1$ follows from Friedman's famous lower bound \cite[Theorem B]{friedman1989analytic} for the usual regulator $R_K$. The second one follows from Siegel's estimate \cite{siegel1969abschatzung}, or \cite[Theorem 1]{louboutin2000explicit}
	$$h_KR_K \leq \dfrac{\omega_K}{2}\left(\dfrac{2}{\pi}\right)^{r_2}\left(\dfrac{e\log |D|}{4(d-1)}\right)^{d-1}\sqrt{|D|}\prod\limits_{\substack{v \in S \\ v  \nmid \infty}}\log\mathcal{N}_{K/\mathbb{Q}}(v),$$
	here, we replace $(1/(d-1))^{d-1}$ with $1$ when $d=1$.
\end{proof}

We will use the following lemma. For the convenience of the readers, we prove it here.
\begin{lemma}
	$\omega_K \leq 2d^2$. Moreover, $\omega_K \leq d^2$ if $K$ contains a primitive $n$-th root of unity for some $n > 6$.
	\label{boundomega}
\end{lemma}
\begin{proof}
	It's sufficient to show that $\varphi(n) \geq \sqrt{n}$ for $n \not = 2, 6$.
	
	For $k \geq 1$, set $f_k(x) := x^k - x^{k-1} - x^{k/2}$, $g_k(x) := x^k - x^{k-1} - \sqrt{2}x^{k/2}$.
	Then 
	$$f_k(x) = x^{(k-1)/2}(x^{(k-1)/2}(x-1) - x^{1/2}) \geq x-1 - x^{1/2} > 0,$$
	if $x \geq 3$. Similarly, $g_k(x) > 0$ if $x \geq 5$ or $k \geq 2, x \geq 3$. 
	
	Let $n = 2^m \prod\limits_{p} p^{e_p}$, where $p$ runs through  all odd prime numbers. If $m = 0$, then 
	$$\varphi(n) = \prod\limits_{e_p \geq 1} (p^{e_p} - p^{e_p-1}) \geq \prod\limits_{e_p \geq 1} p^{e_p/2} = \sqrt{n}.$$ It is similar for the case where $m \geq 2$. 
	
	If $m = 1$, then there exists a prime $q$ such that $q \geq 5, e_q \geq 1$ or $q = 3, e_q \geq 2$. Hence
	\[
	\varphi(n) = \prod\limits_{e_p \geq 1} (p^{e_p} - p^{e_p-1}) \geq \sqrt{2}q^{e_q/2}\prod\limits_{\substack{p \not = q \\e_p \geq 1}} p^{e_p/2} = \sqrt{n}.
	\qedhere
	\]
\end{proof}

\section{Baker's inequality}

In this section, we state Baker's inequality in an explicit form.

\begin{theorem}[Baker's inequality]
	Let $n$ be an integer not less than $2$, $K$ be a number field of degree $d$, $\alpha_1,\cdots, \alpha_n \in K^*$, and $b_1,\cdots, b_n \in \mathbb{Z}$ such that $\alpha_1^{b_1}\cdots\alpha_n^{b_n}\neq 1$. We define $A_1,\cdots, A_n, B_0$ by 
	$$\log A_i:=\max\{\height(\alpha_i),1/d\}, \quad 1 \leq i \leq n;$$
	$$B_0 := \max\{3,|b_1|,\cdots,|b_n|\}.$$ 
	Then for any $v\in M_K$, we have
	\begin{equation}
	|\alpha_1^{b_1}\cdots\alpha_n^{b_n}-1|_v \geq \exp\{-\Upsilon\log A_1\cdots \log A_n\log B_0\},
	\label{baker's main inequality}
	\end{equation}
	where 
	\begin{equation}
	\Upsilon=\begin{cases}
	2^{8n+29} d^{n+2}\log(ed)&\mbox{if $ v|\infty$,}\\
	2^{10n+10}\cdot e^{2n+2}d^{3n+3}p_v^d &\mbox{if $v|p_v<\infty$.}
	\end{cases}
	\end{equation}
	\label{Baker ineq}
\end{theorem}

The proof of this theorem is based on \cite[Corollary 2.3]{matveev2000explicit} and \cite[Main Theorem, page 190-191]{yu2007p}.
For the convenience of readers, we state their results here.

As convention, for a nonzero element $z \in \C$ we set 
$$\log z = \log|z| + \sqrt{-1}\arg z,$$
where $-\pi<\arg z \leq \pi$ is the principal argument of $z$.

\begin{theorem}[{\cite[Corollary~2.3]{matveev2000explicit}}]
	Let $n \in \N^+$, let $K$ be a number field of degree $d$, and let $\alpha_1,\cdots, \alpha_n \in K^*$. Let  $b_1,\cdots, b_n \in \Z$ be such that $\Lambda: = b_1\log\alpha_1 + \cdots + b_n\log\alpha_n \not = 0$. We define $A^*_1,\cdots,A^*_n, B$ by
	$$\log A^*_i = \max\{\height(\alpha_i), \frac{|\log\alpha_i|}{d}, \frac{0.16}{d} \},\quad 1\leq j \leq n,$$
	$$B = \max\{3, \frac{|b_j|\log A^*_j}{\log A^*_n}: 1\leq j \leq n\}.$$
	Then 
	$$\log|\Lambda| \geq -C(n,\varkappa)d^{n+2}\log(ed)\log A^*_1\cdots\log A^*_n\log(eB),$$
	where $C(n,\varkappa) = \min\{\frac{1}{\varkappa}(\frac{1}{2}en)^{\varkappa}30^{n+3}n^{3.5},2^{6n+20}\},$
	\begin{equation*}
	\varkappa =
	\begin{cases}
	1 & \mbox{if $\alpha_1,\cdots , \alpha_n \in \R$,}\\
	2 & \mbox{otherwise.}
	\end{cases}
	\end{equation*}
	\label{matveev's theorem}
\end{theorem}	

\begin{theorem}[\cite{yu2007p} consequence of Main Theorem]
	Keep the notation of Theorem~\ref{matveev's theorem}. We define $A_1,\cdots,A_n, B_0$ by
	$$\log A_i = \max\{\height(\alpha_i), \frac{1}{16e^2d^2}\}, \quad 1\leq i \leq n,$$
	$$B_0 = \max\{3,|b_1|,\cdots,|b_n|\}.$$
	Then for any prime number $p$, and any prime ideal $\mathfrak{p}$ over $p$ in the ring of integers of  $\Q(\alpha_1, \cdots, \alpha_n)$, we have 
	$$\ord_{\mathfrak{p}}(\alpha_1^{b_1}\cdots \alpha_n^{b_n} - 1) < C_0(n,d,\mathfrak{p})\log A_1\cdots \log A_n\log B_0,$$
	where $C_0(n,d,\mathfrak{p}) = (16ed)^{2(n+1)}n^{5/2}\log(2nd)\log(2d)\cdot e^n_{\mathfrak{p}}\frac{p^{f_\mathfrak{p}}}{(f_\mathfrak{p}\log p)^2}$, and $e_{\mathfrak{p}}, f_{\mathfrak{p}}$ are the ramification index and the residue degree at $\mathfrak{p}$ respectively.
	\label{yu's theorem}
\end{theorem}

Now we prove Theorem~\ref{Baker ineq}. The idea comes from \cite[Section 9.4.4]{waldschmidt2013diophantine}.

\begin{proof}[Proof of Theorem~\ref{Baker ineq}]
	If $v|p_v$ for some prime $p_v$, then from Theorem \ref{yu's theorem}, we have 
	$$|\alpha_1^{b_1}\cdots \alpha_n^{b_n} - 1|_v > \exp\{-C_1(n,d,\mathfrak{p})\log A_1\cdots \log A_n\log B_0\},$$
	where
	\[
	C_1(n,d,\mathfrak{p}) = (\frac{\log p_v}{e_\mathfrak{p}})C_0(n,d,\mathfrak{p}) = (16ed)^{2(n+1)}n^{5/2}\log(2nd)\log(2d)\cdot e^{n-1}_{\mathfrak{p}}\frac{p_v^{f_\mathfrak{p}}}{f^2_\mathfrak{p}\log p_v}.
	\]
	We have 
	\begin{align*}
	C_1(n,d,\mathfrak{p})& \leq (16e)^{2(n+1)}d^{2n+2} n^{5/2}\cdot 2nd\cdot 2d\cdot d^{n-1} \cdot p_v^d\\
	& \leq 2^{10n+10}\cdot e^{2n+2}d^{3n+3}p_v^d,
	\end{align*}
	since $n^{7/2} \leq 4^n$.
	
	If $v| \infty$, it is sufficient to bound $|\alpha_1^{b_1}\cdots \alpha_n^{b_n} -1|$. When $|z|\leq 1/2$, the function $\frac{\log(1+z)}{z}$ is holomorphic, then by the maximal modulus principle,  there exists $z_0$ with $|z_0| = 1/2$ such that 
	$$\left|\frac{\log(1+z)}{z}\right| \leq 2|\log(1+z_0)| \leq 2\log2,$$
	where we use the inequality $|\log(1+z_0)| = |\sum\limits_{n=1}^{\infty}\frac{(-1)^{n-1}}{n}z_0^n| \leq \sum\limits_{n=1}^{\infty}\frac{1}{n}|z_0|^n = \log2.$ Hence, for $|z|\leq 1/2$,
	\begin{equation}
	|\log(1+z)|\leq 2\log2|z| \leq 2|z|.
	\label{bound log(1+z)}
	\end{equation}
	To prove  Theorem~\ref{Baker ineq}, without loss of generality, we may assume that $b_i \not = 0$ for $1 \leq i \leq n$, and $A_1 \leq \cdots \leq A_n$, and set $\alpha = \alpha_1^{b_1}\cdots \alpha_n^{b_n} -1$. We need to consider three cases.
	
	(a) If $B_0 \leq 2nd$, with Liouville's inequality, we have 
	$$\height(\alpha) \leq \log 2 + \sum\limits_{i=1}^{n}|b_i|\height(\alpha_i),$$
	$$\log|\alpha| \geq -d\height(\alpha) \geq -d(\log 2 + nB_0\log A_n).$$
	Hence,
	$$|\alpha| \geq \exp\{-(d\log 2+ 2n^2d^2\log A_n)\}.$$
	Since  $1 \leq  d\log A_i$ for $1\leq i \leq n$, and 
	$\log 2 + 2n^2 \leq 2^{8n+29} \log(ed),$ we have 
	$$d\log 2+ 2n^2d^2\log A_n \leq (\log 2 + 2n^2) d^2\log A_n\leq  \Upsilon\log A_1\cdots \log A_n\log B_0.$$
	Hence we have inequality \eqref{baker's main inequality}.
	
	(b) If $B_0 > 2nd$, and $|\alpha| > 1/2$, since $\log 2 \leq 2^{8n+29} \log(ed)$, it is easy to deduce inequality \eqref{baker's main inequality} from this.
	
	(c) If $B_0 > 2nd$, and $|\alpha| \leq 1/2$, this is the main part of the proof. By \eqref{bound log(1+z)}, we have 
	$$|\alpha| \geq \frac{1}{2}|\log(1+\alpha)| = \frac{1}{2}|\log(\alpha_1^{b_1}\cdots\alpha_n^{b_n})| = \frac{1}{2}|\Lambda|,$$
	where $\Lambda = b_0\log(-1) + b_1\log\alpha_1+ \cdots + b_n \log\alpha_n$, $b_0 = 2k$ for some integer $k$. Hence, it is sufficient to bound $|\Lambda|$. 
	
	To use Theorem~\ref{matveev's theorem}, for $1 \leq i \leq n$ we set 
	$$\log A_i^* = \sqrt{\pi^2+1}\cdot \log A_i,$$
	$$\log A_0^* = \frac{\pi}{d},$$
	$$B = B_0^2.$$ 
	We will show that for $1 \leq i \leq n$,  we have
	$$\log A_i^* \geq \max\{\height(\alpha_i),\frac{|\log\alpha_i|}{d}, \frac{0.16}{d}\},$$
	$$\log A_0^* \geq \max\{\height(-1),\frac{|\log(-1)|}{d}, \frac{0.16}{d}\} = \frac{\pi}{d},$$
	$$B \geq \max\{3, \frac{|b_j|\log |A^*_j|}{\log A_n^*}: 0 \leq j \leq n\}.$$ 
	Indeed, notice that for $1\leq i \leq n$, $\log A_i^* \geq\frac{0.16}{d}$, and we have 
	$$|\log \alpha_i|^2 \leq \pi^2 + (\log|\alpha_i|)^2,$$
	$$\frac{\log|\alpha_i|}{d} \leq \height(\alpha_i) \leq \log A_i < \log A^*_i,$$
	so
	$$|\log\alpha_i| \leq (\pi^2 + d^2(\log A_i)^2)^{1/2} \leq \sqrt{\pi^2+1}\cdot d\log A_i.$$
	
	For $\log A_0^*$, it's obvious.
	For $B$, obviously $B \geq 3$. Before showing that $B \geq \frac{|b_j|\log |A^*_j|}{\log A_n^*}$ for $0 \leq j \leq n,$ we bound $b_0$ first. Since $|\alpha| \leq 1/2$, so $|\Lambda| \leq 1$ and
	\begin{align*}
	\pi|b_0| &\leq |\Lambda| + |b_1\log\alpha_1+ \cdots + b_n \log\alpha_n|\\
	&\leq 1 + nB_0\sqrt{\pi^2+1}d\log A_n\\
	&\leq 2\pi ndB_0\log A_n,
	\end{align*}
	here we use the fact that $\sqrt{\pi^2+1} \leq \pi +1$, $1 \leq (\pi -1 )ndB_0\log A_n$.
    Since $B_0 > 2nd \geq 2n$, we have $B=B_0^2 > 2nB_0$,
	$$\frac{|b_0|\log A_0^*}{\log A_n^*} = \frac{\pi|b_0|}{\sqrt{\pi^2+1}\cdot d\log A_n} \leq \frac{2\pi}{\sqrt{\pi^2 + 1}}nB_0 < 2nB_0 < B,$$
	$$\frac{|b_i|\log A_i^*}{\log A_n^*} = \frac{|b_i|\log A_i}{\log A_n} \leq |b_i| \leq B_0 < B$$
	for $1\leq i \leq n$.
	
	By Theorem~\ref{matveev's theorem}, we have
	\begin{align*}
	\log|\Lambda|
	{}&{}\geq -C(n+1,\varkappa)d^{n+3}\log(ed)\log A^*_0\log A^*_1\cdots\log A^*_n\log(eB)\\
	& \geq -3\pi (\pi^2+1)^{n/2}C(n+1,\varkappa)d^{n+2}\log(ed)\cdot \log A_1\cdots\log A_n\log B_0,
	\end{align*}
	and
	{\small
	\begin{align*}
	|\alpha|
	{}&{} \geq \frac{1}{2}|\Lambda| \\
	& \geq \exp\{-(3\pi (\pi^2+1)^{n/2}C(n+1,\varkappa) + \log 2)d^{n+2}\log(ed)\cdot \log A_1\cdots\log A_n\log B_0\}.
	\end{align*}
	}
    Hence, it is sufficient to show that 
	$$3\pi (\pi^2+1)^{n/2}C(n+1,\varkappa) + \log 2 \leq 2^{2n+3}C(n+1,\varkappa) \leq 2^{8n+ 29}.$$
	Indeed, 
	\begin{align*}
	2(\pi^2+1)^{n/2} \left(4\cdot\left(\frac{4}{\sqrt{\pi^2+1}}\right)^n - \frac{3}{2}\pi\right) {}&{} C(n+1,\varkappa) \\
	&\geq  2(\pi^2+1)^{1/2}\left(\frac{16}{\sqrt{\pi^2+1}} - \frac{3}{2}\pi\right)C(2,\varkappa)\\
	&\geq 0.92 \cdot C(2,\varkappa)\\
	& \geq \log 2,
	\end{align*}
	since $0.92\cdot C(2,\varkappa) \geq 0.92\cdot\min\{2^{2.5}e\cdot 30^5, 2^{32}\} \geq \log2$. 
\end{proof}

The following lemma will be used when we apply Theorem~\ref{Baker ineq}.
\begin{lemma}[{\cite[Lemma~2.2]{petho1986products}}]
	Let $b \geq 0, h \geq 1, a> (e^2/h)^h$, and let $x \in \mathbb{R}^+$ be such that 
	$$x-a(\log x)^h - b \leq 0,$$
	then $x < 2^h(b^{1/h}+a^{1/h}\log(h^ha))^h$. In particular, if $h=1$, then $x < 2(b+a\log a).$
	\label{x-alogx-b} 
\end{lemma}

\section{Proof of Theorem~\ref{constant}}

We only consider the case of mixed level, i.e Theorem~\ref{bound}, since if $N$ is a power of some prime $p$, we can replace $N$ by $3N$ if $p =2$, and by $2N$ if $p \not =2$. From the assumption, we have that $N\geq 6$.  

We consider the case where $\Q(\zeta_N) \subset K$ at first, then consider the general case.

For $P \in X_\Gamma(\mathcal{O}_S,j)$, since $j(P) \in \mathcal{O}_S$, we have 
$$\height(j(P))=d^{-1}\sum\limits_{v\in S}d_v\log^+|j(P)|_v \leq \sum\limits_{v\in S}\log^+|j(P)|_v \leq s\log|j(P)|_w,$$
for some $w \in S$. Hence, it suffices to bound $\log|j(P)|_w$. 

If $|j(P)|_w \leq 3500$, then $\height(j(P))\leq 16s$, which is a better bound than that given in Theorem~\ref{bound} when $C = 2^{14}$.

If $|j(P)|_w > 3500$, then by \cite[Proposition 3.3]{sha2014bounding} or \cite[Proposition 3.1]{bilu2011runge}, we have $P\in \Omega_{c,w}$ for some cusp $c$, and $|j(P)|_w \leq 2|q_w(P)^{-1}|_w$ , where  $\Omega_{c,w}$ and $q_w$ are defined in \cite[Section~3]{bilu2011runge}. Hence, we only need to bound $\log|q_w(P)^{-1}|_w$.

Notice that, if $|q_w(P)|_w > 10^{-N}$, then $\log |j(P)|_w \leq 2N\log 10$ and $\height(j(P))< 6sN$, which is better than that given in  Theorem~\ref{bound} when $C=2^{14}$.

In the sequel, we consider the case where $P\in \Omega_{c,w}$ and $|q_w(P)|_w \leq 10^{-N}$.

By the statements in \cite[Page 4507-4508]{sha2014bounding},
there exists a modular unit $W$ on $X_\Gamma$ which is integral over $\mathbb{Z}[j]$, and a constant $\gamma_w \in \mathbb{Q}(\zeta_N)$ such that
$$|\gamma_w^{-1}W(P) - 1|_w \leq 4^{24N^7}|q_w(P)|_w^{1/N},$$
$$\height(\gamma_w)\leq 24N^7\log2,$$ 
and $W(P)$ is a unit of $\OO_{S}$.
Hence $W(P) = \omega\eta_1^{b_1}\cdots\eta_r^{b_r}$ for some $b_1,\cdots, b_r\in \mathbb{Z}$, where $\omega$ is a root of unity and $\{\eta_1,\cdots,\eta_r\}$ is a fundamental system of $S$-units from \cite[Proposition~4.1]{sha2014bounding}. We set
$$\Lambda= \gamma_w^{-1}W(P) = \eta_0\eta_1^{b_1}\cdots\eta_r^{b_r},$$
where $\eta_0 = \omega\gamma_w^{-1}$. Then we have
\begin{equation}\label{lambdaupbound}
|\Lambda-1|_w\leq 4^{24N^7}|q_w(P)|_w^{1/N}.
\end{equation}
If $\Lambda\not = 1$, we will use this upper bound and the lower bound from Theorem~\ref{Baker ineq} to get a bound of $|q_w(P)|_w$ which gives an upper bound of $\height(j(P))$.
For the case where $\Lambda = 1$, see \cite[Section 8]{sha2014bounding}. 

To state the following lemma, we set $r^{r} = 1$ when $r = 0$, i.e $s = 1$.

\begin{lemma}
	If $\Q(\zeta_N) \subset K$ and $\Lambda\neq 1$, then we have
	$$\height(j(P)) \leq 40dsr^{2r}\zeta^rN^8\tilde{\Upsilon} R(S)\log(d^2sr^{4r}\zeta^{s}N^{16}\tilde{\Upsilon}R(S)),$$
	where $\tilde{\Upsilon} = 2^{13s+22}d^{2s+3}\ell^d$,  and $\zeta$ has been defined in Section~\ref{S-regulator}.
	\label{bound1}
\end{lemma}
\begin{proof}
	
	We define $A_0,\cdots, A_r, B_0$ by 
	$$\log A_i:=\mbox{max}\{\mbox{h}(\eta_i),1/d\}, 0 \leq i \leq r;$$
	$$B_0 := \mbox{max}\{3,|b_1|,\cdots,|b_r|\}.$$ 
	Since $\Lambda= \eta_0\eta_1^{b_1}\cdots\eta_r^{b_r} \not= 1$, by Theorem~\ref{Baker ineq}, we have
	$$|\Lambda-1|_w \geq \exp\{-{\Upsilon}\log A_0\cdots \log A_r\log B_0\},$$
	where 
	\begin{equation}
	{\Upsilon}=\begin{cases}
	2^{8s+29} d^{s+2}\log(ed),&\mbox{if $ w|\infty$,}\\
	2^{10s+10}\cdot e^{2s+2}d^{3s+3}p_w^d, &\mbox{if $w|p_w<\infty$.}
	\end{cases}
	\end{equation}
	Obviously $2^{10s+19}\cdot 2^{3s+3}d^{3s+3}\ell^d = 2^{13s+22}d^{3s+3}\ell^d$ is larger than ${\Upsilon}$ in each one of the cases since $d\geq 2, s \geq 1$, so we can take ${\Upsilon} = 2^{13s+22}d^{3s+3}\ell^d$. 
	
	By \eqref{lambdaupbound}, we have 
	$$\exp\{-{\Upsilon}\log A_0\cdots \log A_r\log B_0\} \leq 4^{24N^7}|q_w(P)|_w^{1/N},$$
    that is
	\begin{equation}
	\log|q_w(P)^{-1}|_w\leq N{\Upsilon}\log A_0\cdots \log A_r\log B_0 + 48N^8\log 2.
	\label{ineqqw}
	\end{equation}
	By \cite[Proposition 4.1]{sha2014bounding}, we have $\zeta \height(\eta_k) \geq 1/d$ and $\zeta \geq 1$, so 
	$$\log A_k \leq \zeta \height(\eta_k), \ \ k=1,\cdots,r,$$
	$$\log A_1 \cdots \log A_r \leq d^{-r}r^{2r}\zeta^rR(S).$$
	Notice that the both sides are $1$ when $r = 0$. 
	On the other hand, since 
	$$\height(\eta_0) = \height(\gamma_w) \leq 24N^7\log 2,$$
	we have $$\log A_0 \leq 24N^7\log2.$$
	
	For $B_0$, we set $B^* = \mbox{max}\{|b_1|,\cdots,|b_r|\}$ if $r \geq$ 1, and $B^* = 0$ if $r = 0$. By \cite[Corollary~4.2 and Proposition~6.1]{sha2014bounding}  we have  
	\begin{equation}
	\begin{aligned}
	B^*&\leq 2dr^{2r}\zeta \mbox{h}(W(P))\\
	&\leq 2dr^{2r}\zeta(2sN^8\log|q_w^{-1}(P)|_w + 94sN^8\log N),
	\end{aligned}
	\label{B0}
	\end{equation}
    so $$B_0 \leq 2dr^{2r}\zeta(2sN^8\log|q_w^{-1}(P)|_w + 94sN^8\log N).$$
	We write 
	$$\alpha = 4dsr^{2r}\zeta N^8,$$
	$$\beta = 188dsr^{2r}\zeta N^8\log N = 47\alpha\log N,$$
	$$C_1 = \alpha N{\Upsilon}\log A_0\cdots \log A_r,$$
	$$C_2 = 48\alpha N^8\log 2 + \beta.$$
	Hence, inequalities \eqref{ineqqw} and \eqref{B0} yield 
	$$\alpha \log|q_w(P)^{-1}|_w + \beta \leq C_1\log (\alpha \log|q_w(P)^{-1}|_w + \beta) + C_2.$$
	By Lemma~\ref{x-alogx-b}, we obtain
	$$\alpha \log|q_w(P)^{-1}|_w + \beta \leq 2(C_1\log C_1 + C_2).$$
	Hence,
	$$\log|q_w(P)^{-1}|_w \leq 2\alpha^{-1}C_1\log C_1 + \alpha^{-1}(2C_2-\beta),$$
	$$\log|j(P)|_w \leq \log2|q_w(P)^{-1}|_w \leq 2\alpha^{-1}C_1\log C_1 + \alpha^{-1}(2C_2-\beta) + \log 2,$$
	so we have 
	$$\height(j(P)) \leq 2s\alpha^{-1}C_1\log C_1 + s\alpha^{-1}(2C_2-\beta) + s\log 2.$$
	Next we bound each term on the right-hand side:
	\begin{align*}
	2s\alpha^{-1}C_1 {}&{} \log C_1 \\
	& = 2sN{\Upsilon}\log A_0\cdots\log A_r\log(4dsr^{2r}\zeta N^9{\Upsilon}\log A_0\cdots\log A_r)\\
	& \leq 48\log 2\cdot d^{-r}sr^{2r}\zeta^rN^8{\Upsilon} R(S)\log (96\log 2\cdot d^{-r+1}sr^{4r}\zeta^{r+1}N^{16}{\Upsilon} R(S))\\
	&\leq 39d^{-r}sr^{2r}\zeta^rN^8{\Upsilon} R(S)\log (d^{-r+1}sr^{4r}\zeta^{r+1}N^{16}{\Upsilon} R(S)),
	\end{align*}
	here we use the fact that $48\log2 \times \log(96\log2) \leq 140 < 5\log (d^{-r+1}{\Upsilon})$; we also have 
	\begin{equation*}
	\begin{aligned}
	s\alpha^{-1}(2C_2-\beta)+s\log 2& = 96\log 2\cdot sN^8 + 47s\log N +s\log 2\\
	& \leq 98\log 2\cdot sN^8.\\
	\end{aligned}
	\end{equation*}
	After replacing $d^{-s}{\Upsilon} = 2^{13s+22}d^{2s+3}\ell^d$ by $\tilde{\Upsilon}$, we have
	\[
	\height(j(P)) \leq 40dsr^{2r}\zeta^rN^8\tilde{\Upsilon}R(S)\log (d^2sr^{4r}\zeta^sN^{16}\tilde{\Upsilon}R(S)).
	\qedhere
	\]
\end{proof}

We will use the bound $\zeta \leq 2^{13}(\log d)^3$ subsequently. If $d =2$, $$\zeta = \dfrac{(\log 6)^3}{2} = \dfrac{(\log_2 6)^3}{2} (\log d)^3 \leq 2^4(\log d)^3;$$
if $d \geq 3$, then
\[
\zeta
= 4\left(\dfrac{\log d}{\log\log d}\right)^3
\leq 4\left(\dfrac{\log d}{\log\log 3}\right)^3
\leq 4809(\log d)^3
\leq 2^{13} (\log d)^3.
\]

By Lemma~\ref{bound RS} and Lemma~\ref{boundomega}, we have
\begin{gather*}
R(S)\leq \dfrac{\omega_K}{2}\dfrac{1}{(d-1)^{d-1}}\left(\log |D|\right)^{d-1}\sqrt{|D|}\prod\limits_{\substack{v \in S \\ v  \nmid \infty}}\log\mathcal{N}_{K/\mathbb{Q}}(v), \\
\omega_K \leq 2d^2, \\
\log R(S)
\leq \log(\dfrac{\omega_K}{2}) + d\log|D|+s\log (d\ell)
\leq 2\log d + d\log|D|+s\log(d\ell).
\end{gather*}
We have $d \leq 2s$ and $\log s \leq s/2$. Then we have
\begin{align*}
\log \tilde{\Upsilon} &= (13s+22)\log 2 + (2s+3)\log d + d\log \ell \\
&\leq (15s+25)\log 2 + (2s+3)\log s + d\log \ell \\
&\leq 28s + (s+2)s + s\ell \\
& \leq 32s^2\ell
\end{align*}
and
\begin{align*}
\log (d^2sr^{4r}\zeta^{r+1}N^{16} {}&{} \tilde{\Upsilon} R(S)) \\
\leq & 2\log d+ 4s\log s+ 13s\log 2 + 3s\log\log d +16\log N \\ 
& +\log\tilde{\Upsilon} + 2\log d + d\log|D| +s\log(d\ell)\\
\leq & 2s+ 2s^2+ 10s + 2s^2 +16\log N + 32s^2\ell + 2s  + 2s\log|D| +s^2\ell\\
\leq & 8N + 2s\log|D| + 51s^2\ell\\
\leq & 61s^2\ell N \log|D|\\
\leq & 2^{6}s^2N\ell\log|D|.
\end{align*}
Hence combining with Lemma~\ref{bound1}, we have
{\small
\begin{align} \label{bound contain root}
\height(j(P))& \notag \\
\leq &{\ } 2^6\cdot ds^{2s-1}\zeta^rN^8\tilde{\Upsilon} R(S)\log(d^2sr^{4r}\zeta^{r+1}N^{16}\tilde{\Upsilon}R(S)) \notag \\
\leq &{\ } 2^{26s+15}\cdot d^{2s+4}(\log d)^{3r}s^{2s-1}N^8\ell^d\dfrac{\omega_K}{2}\dfrac{1}{(d-1)^{d-1}}\left(\log |D|\right)^{d-1}\sqrt{|D|} \notag \\
&\cdot \prod\limits_{\substack{v \in S \\ v  \nmid \infty}}\log\mathcal{N}_{K/\mathbb{Q}}(v) \cdot(2^{6}s^2N\ell\log|D|) \notag \\
= &{\ } 2^{26s+20}d^{2s+4}(\log d)^{3r}s^{2s+1}N^9\ell^{d+1}\omega_K\dfrac{1}{(d-1)^{d-1}}(\log|D|)^d\sqrt{|D|}\prod\limits_{\substack{v \in S \\ v  \nmid \infty}}\log\mathcal{N}_{K/\mathbb{Q}}(v).
\end{align}
}

Next we deal with the general case.
Set $\tilK = K\cdot\mathbb{Q}(\zeta_N) = K(\zeta_N)$. Let $\tilS$ be the set consisting of the extensions of the places from $S$ to $\tilK$, that is,
$$\tilS= \{\tilv\in M_\tilK: \tilv|v, v \in S\}.$$
Then $ P \in X_\Gamma(\mathcal{O}_{\tilS},j)$. Put $\tild = [\tilK:\mathbb{Q}]$, $\tils = |\tilS|$,  $\tilde{r} = \tils-1$,  
and let $\tilD$ be the absolute discriminant of $\tilK$.
\begin{lemma}
\begin{gather*}
N-\varphi(N)\geq 4, \\
\tils\leq s\varphi(N), \\
\tild\leq d\varphi(N), \\
\omega_\tilK \leq 2d^2\varphi(N)^2, \\
|\tilD|\leq N^{dN}|D|^{\varphi(N)}, \\
\prod\limits_{\substack{v \in \tilS \\ v  \nmid \infty}}\log\mathcal{N}_{\tilK/\mathbb{Q}}(v) \leq 4^{s\varphi(N)}\left(\prod\limits_{\substack{v \in S \\ v  \nmid \infty}}\log\mathcal{N}_{K/\mathbb{Q}}(v)\right)^{\varphi(N)}.
\end{gather*}
\end{lemma}
\begin{proof}
	The first three inequalities come directly from the definition of $\tilK$ and $\tilS$ and $N \geq 6$  has at least two prime factors. The fourth inequality comes from $\omega_\tilK \leq 2\tild^2 \leq 2d^2\varphi(N)^2$.

	Let $D_{\tilK/K}$ be the relative discriminant of $\tilK/K$. We have 
	$$\tilD=\mathcal{N}_{K/\mathbb{Q}}(D_{\tilK/K})D^{[\tilK:K]}.$$
	We denote by $\mathcal{O}_{K}$ and $\mathcal{O}_\tilK$ the ring of integers of $K$ and $\tilK$, respectively. Since $\tilK=K(\zeta_N)$, we have 
	$$\mathcal{O}_{K} \subset \mathcal{O}_{K}(\zeta_N) \subset \mathcal{O}_\tilK.$$
	Note that the absolute value of the discriminant of the polynomial $x^N - 1$ is $N^N$, we obtain 
	$$D_{\tilK/K}|N^N,$$
	so
	$$|\mathcal{N}_{K/\mathbb{Q}}(D_{\tilK/K})| \leq N^{dN}.$$
	Hence,
	$$|\tilD|\leq N^{dN}|D|^{\varphi(N)}.$$

	Notice that $\tilK/K$ is Galois. Let $v$ be a non-Archimedean place of $K$, and let $v_1,\dots, v_g$ be all its extensions to $\tilK$ with residue degree $f$ over $K$. Then $gf\leq [\tilK:K] \leq \varphi(N)$, which implies $g\log_2 f\leq gf\leq \varphi(N)$, i.e. $f^g \leq 2^{\varphi(N)}$. Note that $2\log\mathcal{N}_{K/\mathbb{Q}}(v) > 1$ and  $\mathcal{N}_{\tilK/\mathbb{Q}}(v_k) = \mathcal{N}_{K/\mathbb{Q}}(v)^{f}$ for $1 \leq k \leq g$, $g\leq \varphi(N)$, we have 
	
	\begin{equation*}
	\begin{aligned}
	\prod\limits_{k=1}^g\log\mathcal{N}_{\tilK/\mathbb{Q}}(v_k) & \leq 2^{\varphi(N)}(\log\mathcal{N}_{K/\mathbb{Q}}(v))^g \\
	& \leq 2^{\varphi(N)}(2\log\mathcal{N}_{K/\mathbb{Q}}(v))^g\\
	& \leq 4^{\varphi(N)}(\log\mathcal{N}_{K/\mathbb{Q}}(v))^{\varphi(N)}.
	\end{aligned}
	\end{equation*}
	Hence
	\[
	\prod\limits_{\substack{v \in \tilS\\  v\nmid \infty}}\log\mathcal{N}_{\tilK/\mathbb{Q}}(v) \leq 4^{s\varphi(N)}\left(\prod\limits_{\substack{v \in S\\  v\nmid \infty}}\log\mathcal{N}_{K/\mathbb{Q}}(v)\right)^{\varphi(N)}.
	\qedhere
	\]
\end{proof}

Combine the lemma above with the bound \eqref{bound contain root}, we have
{\small
\begin{align*}
\height(j&(P)) \\
&\leq 2^{26\tils+20}\tild^{2\tils+4}(\log \tild)^{3\tilde{r}}\tils^{2\tils+1}N^9\ell^{\tild+1}\omega_\tilK\dfrac{1}{(\tild-1)^{\tild-1}}(\log|\tilD|)^\tild\sqrt{|\tilD|}\prod\limits_{\substack{v \in \tilS \\ v  \nmid \infty}}\log\mathcal{N}_{\tilK/\mathbb{Q}}(v)\\
&\leq 2^{28s\varphi(N)+21}d^{2s\varphi(N)+6}(\log d\varphi(N))^{3s\varphi(N)}s^{2s\varphi(N)+1}\varphi(N)^{4s\varphi(N)+7}N^{9}\ell^{d\varphi(N)+1}\Delta(N)\\
&\leq 2^{28sN}d^{2sN}(\log dN)^{3sN}s^{2sN}N^{4sN}\ell^{dN}\Delta(N)\\
&\leq (2^{14}dsN^2)^{2sN}(\log dN)^{3sN}\ell^{dN}\Delta(N).
\end{align*}
}
This completets the proof of Theorem~\ref{constant}.

\section*{Acknowledgements}
The research of the author is supported by the China Scholarship Council. The author also thanks his supervisors Yuri Bilu and Qing Liu for helpful discussions and valuable suggestions about this paper, and careful reading.

%%% References

\EditInfo{%
    December 27, 2019}{%
    April 04, 2020}{%
    Yuri Bilu}

\end{paper}